\documentclass[12pt]{article}
\usepackage{amsthm,amssymb}
\def\del#1{}
\def\dell#1{}

\newcommand{\limfunc}[1]{{\rm \mathop{#1}}}
\newcommand{\text}[1]{{\mbox {#1}}}
\newcommand{\func}[1]{{\rm #1} \,}
\newtheorem{theorem}{Theorem}
\newtheorem{lemma}{Lemma}
\newtheorem{proposition}{Proposition}
\newtheorem{remark}{Remark}
\newtheorem{corollary}{Corollary}

\newlength{\rig}
\newlength{\rigg}
\newlength{\hei}

\newcommand{\dokaend}{\hfill$\square$ \vskip6truept}

\begin{document}
\setcounter{page}{1}
\author{{\bf Anton Savin $^*$
}\\
Moscow State University\\
e-mail: antonsavin@mtu-net.ru\\
{\bf Bert-Wolfgang Schulze $^*$}\\
Potsdam University\\
e-mail: schulze@math.uni-potsdam.de\\
{\bf Boris Sternin
\thanks{
Supported by DAAD via IQN ``Kopplungsprozesse und ihre Strukturen in der Geo- und Biosph\"are'',
by DFG via a project with the Arbeitsgruppe ``Partielle Differentialgleichungen und Komplexe
Analysis,''
Universit\"at Potsdam, by RFBR grants Nos.~02-01-00118,
02-01-06515 and~00-01-00161, and by the European Training Network Program ``Geometric
Analysis''.}
}\\
Independent University of Moscow\\
e-mail: sternine@mtu-net.ru
}
\title{\bf
Elliptic Operators in Subspaces and the Eta Invariant
}

\maketitle

\begin{abstract}
In this paper we obtain a formula for the fractional part of the $\eta$-invariant for elliptic
self-adjoint operators{ in topological terms}. The computation of the $\eta$-invariant is based
on the index theorem for elliptic operators in subspaces obtained in \cite{SaSt1}, \cite{SaSt2}.
We also apply the $K$-theory with coefficients $\Bbb{ Z}_n.$ In particular, it is shown that the
group $K(T^*M,\Bbb{ Z}_n)$ is realized by elliptic operators (symbols)\dell{,} acting in
appropriate subspaces.
\end{abstract}


\vspace{15mm}

{\bf 2000 AMS classification}: Primary 58J28, Secondary 58J20, 58J22,
19K56

\vfill

\newpage


\section*{Introduction}
\addcontentsline{toc}{section}{Introduction}

The spectral $\eta $-invariant of an elliptic self-adjoint operator on a
closed manifold was introduced by Atiyah, Patodi, and Singer
\cite{APS1}. It appeared as a nonlocal contribution to an index formula for
manifolds with boundary obtained via the heat equation method. From the
very moment of its introduction, it was clear that this spectral invariant
in the general case is neither an invariant of the principal symbol of the
operator\dell{,} nor a homotopy invariant of the operator itself. More
precisely, for a one-parameter family $A_t$ of elliptic self-adjoint
operators, the function $\eta \left( A_t\right) $ is {\em piecewise
smooth} with  jumps at those values of $t$, where
some eigenvalue of the operator in the family changes its sign.

P. Gilkey \cite{Gil7} observed that for {\em differential operators}
satisfying the parity condition
\begin{equation}
\limfunc{ord}A+\dim M\equiv 1\left( \func{mod}2\right) ,  \label{conda}
\end{equation}
the $\eta $-invariant of a one-parameter family is
a {\em piecewise constant} function. In particular, in this case the
fractional part of the spectral $\eta $-invariant
is in fact a homotopy invariant\dell{,}
depending on the principal symbol of the operator only.
This rises the problem of computing this invariant
in topological terms and the problem of finding nontrivial
geometric examples.

The situation is rather well understood for even-dimensional manifolds.
In this case the famous first order Dirac type operators satisfy the
parity condition (\ref{conda}). These were studied by P.~Gilkey \cite{Gil8}.
He proved that the $\eta$-invariant takes
dyadic rational values. Nontrivial $\eta$-invariants were computed on
some nonorientable ${\rm pin}^c$ manifolds, e.g. $\!\Bbb{RP}^{\,2n}$.
This fractional invariant is important in topology and differential
geometry (e.g., see \cite{Stol1,BaGi1,Gil10}).

On odd-dimensional manifolds, P.~Gilkey showed in \cite{Gil7} that the
fractional part of the $\eta$-invariant defines a homomorphism
$$
\left.K(P^*M)\right/K(M)\longrightarrow\Bbb{Z}\left[\frac12\right]/\Bbb{Z},
$$
where $P^*M=S^*M/\Bbb{Z}_2$ denotes the projective space bundle of $M$.
Moreover, he introduced a class of second order operators on $M$
that define nontrivial elements of the latter $K$-group  and proposed a
problem of computing their $\eta$-invariants.

The aim of the present paper is to obtain a topological formula for
the fractional part of the $\eta$-invariant for operators satisfying
the parity condition.

The $\eta $-invariant of an operator $A$ satisfying condition (\ref{conda})
is completely determined by the nonnegative spectral subspace
$\widehat{L}_{+}\left(A\right)$ of this operator, while the fractional
part of the invariant is determined by the so-called {\em symbol
of the subspace}. This is a vector bundle on the cospheres $S^{*}M$ over the
manifold\dell{,} generated by the positive eigenspaces of the self-adjoint
symbol $\sigma \left( A\right) $.
{First,} this enables one to {identify}
the $\eta $-invariant of self-adjoint elliptic
operator{s}  {with} a dimension-type functional (see \cite{SaSt1})
{on}  the corresponding (infinite-dimensional) spectral subspace{s}.
Second, we can apply the index formula for
elliptic operators\dell{,} in subspaces \cite{SaSt1,SaSt2}.
The index formula reduces the computation to {the} ''index modulo $n$''
problem for operators  in subspaces. The term ''modulo $n$'' here
expresses the fact that in this case the index of an elliptic operator,
being reduced modulo $n$, becomes an invariant of the
principal symbol of the operator. It turns out that {such} elliptic
operators on a closed manifold\dell{,} define the $K$-theory with
coefficients in $\Bbb{ Z}_n.$ In particular, the index is computed
modulo $n$ by the direct image mapping in $K$-theory.

The fractional part of the $\eta$-invariant was first computed
in the classical Atiyah--Patodi--Singer paper \cite{APS3} for
operators with coefficients in flat bundles. However, their result does
not apply to our operators, since there is no flat bundle available.
Let us also mention that although we do not rely on the results of
\cite{APS3}, there are strong parallels between the two proofs.
For instance,  the index formula in subspaces plays the same role as
the Atiyah--Patodi--Singer index formula for trivialized flat bundles,
while  the mod $n$-index formula of the present paper generalizes the
mod$n$-spectral flow formula.

{\bf Added January 2002.} The main formula for the fractional part of the
$\eta$-invariant was used to find {\em nontrivial} $\eta$-invariants
for some new second order geometric operators in the papers
\cite{SaSt7,SaSt9}. Thus, the formula for the fractional part of the
$\eta$-invariant on odd-dimensional manifolds is not empty and
the problem of finding even-order operators in odd dimensions with
nontrivial $\eta$-invariant is solved.

The results of the paper were reported at the international conferences
''Operator Algebras and Asymptotics on Manifolds with Singularities,''
Warsaw, and "Jean Leray 99", Karlskrona, Sweden. We are thankful
to the participants of Professor A.~S.~Mishchenko
seminar at the Moscow State University for useful discussions.
We would like to thank the referees for their critical remarks
on the preliminary version of the paper.

\section{Subspaces and index formulas}

{\bf 1. }{\bf Pseudodifferential subspaces. The dimension functional.}
Spaces defined by pseudodifferential projections
on a smooth closed manifold $M$
were considered in~\cite{SaSt1,SaSt2}. More precisely,
a subspace
\[
\widehat{L}\subset C^\infty \left( M,E\right)
\]
in the space of smooth sections of a vector bundle $E$ on $M$
is said to be {\em pseudodifferential\/} if it is the range
\[
\widehat{L}=\func{Im}P,\quad P:C^\infty \left( M,E\right) \rightarrow
C^\infty \left( M,E\right)
\]
for some pseudodifferential projection $P$ of order zero.
The principal symbol of the projection defines the vector bundle
\begin{equation}
L=\func{Im}\sigma \left( P\right) \subset \pi ^{*}E\in \limfunc{Vect}\left(
S^{*}M\right).  \label{alal}
\end{equation}
It is the {\em symbol of the subspace}. Here $\pi
:S^{*}M\rightarrow M$ is the projection for the cosphere bundle.

On the cotangent bundle $T^{*}M$, we consider the antipodal involution
\[
\alpha :T^{*}M\longrightarrow T^{*}M,\quad \alpha \left( x,\xi \right)
=\left( x,-\xi \right) .
\]
A subspace $\widehat{L}\subset C^\infty \left( M,E\right) $ is
said to be {\em even (odd)} with respect to $\alpha $
if the symbol $L$ is invariant (antiinvariant):
\begin{equation}
L=\alpha ^{*}L,\quad \text{or\quad }L\oplus \alpha ^{*}L=\pi ^{*}E.
\label{omg}
\end{equation}
We point out that both equalities in this formula are equalities
of subbundles in the ambient bundle $\pi ^{*}E.$
Denote the semigroups of even(odd)
subspaces by $\widehat{\limfunc{Even}}\left( M\right) $ $\left( \widehat{%
\limfunc{Odd}}\left( M\right) \right) .$ The symbols of even (odd) subspaces
will be referred to as even (odd) bundles for brevity.

It turns out that if the parities of the subspaces and of the
dimension of $M$ are opposite, then the subspaces have a homotopy
invariant similar to the dimension for finite-dimensional vector spaces.
More precisely, the following theorem holds.

\begin{theorem}
{\em \cite{SaSt1,SaSt2}}
There is a unique additive functional
\[
d:\widehat{\limfunc{Even}}\left( M^{odd}\right) \rightarrow \Bbb{ Z}\left[
\frac 12\right] ,\quad \text{or \quad }d:\widehat{\limfunc{Odd}}\left(
M^{ev}\right) \rightarrow \Bbb{ Z}\left[ \frac 12\right]
\]
with the following properties:
\end{theorem}

\begin{enumerate}
\item  {\em (invariance) }$d\left( U\widehat{L}\right) =d\left( \widehat{L}%
\right) ${\em \ for all invertible pseudodifferential operators }$U$ {\em %
with even principal symbol: }$\alpha ^{*}\sigma \left( U\right) =\sigma
\left( U\right) ;$

\item  {\em (relative index) }$d\left( \widehat{L}_1\right) -d\left(
\widehat{L}_2\right) =\limfunc{ind}\left( \widehat{L}_1,\widehat{L}_2\right)
${\em \ for two subspaces with coinciding principal symbols;
\footnote{The relative index of subspaces
$\limfunc{ind}\left( \widehat{L}_1,\widehat{L}%
_2\right) $ is expressed via the Fredholm index \cite{BDF1}
\[
\limfunc{ind}\left( \widehat{L}_1,\widehat{L}_2\right) =\limfunc{ind}\left(
P_2:\func{Im}P_1\rightarrow \func{Im}P_2\right),
\]
where the projections $P_{1,2}$ define $\widehat{L}_{1,2}$.}}

\item  {\em (complement) }$d\left( \widehat{L}\right) +d\left( \widehat{L}%
^{\perp }\right) =0,${\em \ \quad where }$\widehat{L}^{\perp }${\em \
denotes the orthogonal complement of} $\widehat{L}.$
\end{enumerate}

\begin{corollary}
The functional $d$ is a homotopy invariant of the subspace, while its
fractional part is an invariant of the symbol of the subspace.
\end{corollary}

Indeed, the homotopy invariance follows from the invariance property.
Moreover, it follows from the relative index property that the fractional
part is determined by the symbol of the subspace.

\noindent {\bf 2. }
{\bf Dimension functional and $\eta$-invariant.}
The functional $d$ can be expressed in terms of the
Atiyah--Patodi--Singer $\eta $-invariant.

Namely, for an elliptic self-adjoint operator
\[
A:C^\infty \left( M,E\right) \longrightarrow C^\infty \left( M,E\right)
\]
of a positive order
consider the subspace $\widehat{L}_{+}\left( A\right) \subset C^\infty
\left( M,E\right) $ generated by the eigenvectors of $A$
corresponding to nonnegative eigenvalues. It is well known
(e.g., see \cite{APS1}) that the spectral
projection $P_{+}\left( A\right) $ on this
subspace is a pseudodifferential operator of order zero.
Thus, the subspace $\widehat{L}_{+}\left( A\right) $ is pseudodifferential.
The symbol $L_{+}\left( A\right)$ of the subspace can be
explicitly calculated:
\[
L_{+}\left( A\right) =\func{Im}\sigma \left( P_{+}\left( A\right) \right)
\subset \pi ^{*}E\in \limfunc{Vect}\left(S^*M\right),
\]
where the principal symbol $\sigma\left(P_{+}\left(A\right)\right)$
of the projection is
equal to the spectral projection for the principal symbol $\sigma \left(
A\right) $:
$\sigma \left( P_{+}\left( A\right) \right) =P_{+}\left( \sigma \left(
A\right) \right).$

Thus, if $A$ is a differential operator, then the subspace $\widehat{L}%
_{+}\left( A\right) $ is either even or odd, according to the parity of
operator's order. The same property holds for a class of pseudodifferential
operators introduced in \cite{Gil7}: these are classical
pseudodifferential operators with homogeneous terms in the
asymptotic expansion of the symbol possessing the $\Bbb{ R}_{*}$-invariance
(cf. \cite{NSS3}):
\begin{equation}
\sigma \left( A\right) \left( x,\xi \right) \sim \sum_{j=o}^\infty
a_{d-j}\left( x,\xi \right) ,\qquad a_k\left( x,-\xi \right) =\left(
-1\right) ^ka_k\left( x,\xi \right) ,\quad \text{for all }k\leq d.
\label{dopus}
\end{equation}
For this class of operators, the functional $d$ is equal to the
$\eta $-invariant.

\begin{theorem}
\label{lab10a}
{\em  \cite{SaSt1,SaSt2}}
For the nonnegative spectral subspace $\widehat{L}_+(A)$ of an elliptic
self-adjoint operator $A$ satisfying {\em (\ref{dopus})} one has
\begin{equation}
d\left( \widehat{L}_+(A)\right) =\eta \left( A\right)  \label{dash}
\end{equation}
provided the order of $A$ and the dimension of the manifold
have opposite parities.
\end{theorem}

According to this result, to compute the fractional part of the
$\eta$-invariant it suffices to compute the fractional part of
the functional $d$  in terms of the symbol of the subspace.
An important ingredient of the computation is the
index formula for elliptic operators in subspaces.

\noindent {\bf 3. }{\bf Elliptic theory in subspaces.}
Let $\widehat{L}_{1,2}\subset C^\infty \left(
M,E_{1,2}\right) $ be two subspaces. Consider a
pseudodifferential operator
\[
D:C^\infty \left( M,E_1\right) \longrightarrow C^\infty \left( M,E_2\right)
\]
in the ambient spaces. If it preserves the subspaces:
$D\widehat{L}_1\subset \widehat{L}_2,$ then the restriction
\begin{equation}
D:\widehat{L}_1\longrightarrow \widehat{L}_2  \label{ttwostar}
\end{equation}
is called an {\em operator acting in subspaces}.
In this case the principal symbol $\sigma \left( D\right) $
restricts to a vector bundle homomorphism
\begin{equation}
\sigma \left( D\right) :L_1\longrightarrow L_2  \label{ssymbl}
\end{equation}
over $S^{*}M$. This is called the {\em symbol of the operator in subspaces}.

It is proved in \cite{ScSS18} that the closure
\[
D:H^s\left( M,E_1\right) \supset \widehat{L}_1\longrightarrow \widehat{L}%
_2\subset H^{s-m}\left( M,E_2\right) ,\qquad m=\limfunc{ord}D,
\]
of (\ref{ttwostar}) in the Sobolev norms
defines a Fredholm operator if  the symbol
(\ref{ssymbl}) is {\em elliptic}, i.e., a vector bundle isomorphism.

For elliptic operators the following index formula was
obtained in \cite{SaSt1,SaSt2}.

\begin{theorem}
One has
\begin{equation}
\limfunc{ind}\left( D,\widehat{L}_1,\widehat{L}_2\right) =\frac 12\limfunc{%
ind}\widetilde{D}+d\left( \widehat{L}_1\right) -d\left( \widehat{L}_2\right)
,  \label{inde}
\end{equation}
where $D:\widehat{L}_1\longrightarrow \widehat{L}_2,$ $\widehat{L}%
_{1,2}\subset C^\infty \left( M,E_{1,2}\right) $ is an elliptic operator in
subspaces of the same parity:
$
\widehat{L}_{1,2}\in \widehat{\limfunc{Even}}\left( M^{odd}\right) \;
or\; \widehat{\limfunc{Odd}}\left( M^{ev}\right),
$
while elliptic operator
\[
\widetilde{D}:C^\infty \left( M,E_1\right) \longrightarrow C^\infty \left(
M,E_2\right)
\]
has principal symbol
\[
\sigma \left( \widetilde{D}\right) =\sigma \left( D\right) \oplus \alpha
^{*}\sigma \left( D\right) :L_1\oplus \alpha ^{*}L_1\longrightarrow
L_2\oplus \alpha ^{*}L_2
\]
for odd subspaces. In the case of even subspaces, $\widetilde{D}$ is
defined as
\begin{eqnarray*}
\widetilde{D}:C^\infty \left( M,E_1\right) \longrightarrow C^\infty \left(
M,E_1\right) ,
\vspace{2mm}
\\
\sigma \left( \widetilde{D}\right) =\left[ \alpha ^{*}\sigma \left( D\right)
\right] ^{-1}\sigma \left( D\right) \oplus 1:L_1\oplus L_1^{\perp
}\longrightarrow L_1\oplus L_1^{\perp }.
\end{eqnarray*}
\end{theorem}

\section{Computation of the $\eta$-invariant\label{secsimp}}

\begin{proposition}
\label{thpow}
For a subspace $\widehat{L}\in \widehat{\limfunc{Even}}\left( M^{odd}\right)
$ \ or $\ \widehat{\limfunc{Odd}}\left( M^{ev}\right) $\ with
symbol $L$, there exists a positive integer $N$ such that
the direct sum $2^NL$ on $S^{*}M$ can be lifted from the base $M.$\
That is, for some
vector bundle $F\in \limfunc{Vect}\left( M\right) $\ there exists an
isomorphism
\begin{equation}
\sigma :2^NL\longrightarrow \pi ^{*}F,\qquad 2^NL=
{\underbrace{L\oplus \ldots \oplus L}_{2^N\text{copies}}}.  \label{omega}
\end{equation}
\end{proposition}

\noindent {\em Proof. }The part of the theorem, pertaining to even
subspaces, follows from \cite{Gil7}, where it is shown that for an
odd-dimensional $M$ the projective space bundle
\[
P^{*}M=\left. S^{*}M\right/ \left\{ \left( x,\xi \right) \sim \left( x,-\xi
\right) \right\}
\]
has the same $K$-theory groups as $M$, except for the 2-torsion.
The isomorphism, modulo 2-torsion, is established by the natural projection
$\pi _P:P^{*}M\rightarrow M.$
More precisely, $\ker \pi _P^{*}=0,$ and $\limfunc{coker}\pi _P^{*}$ is a
2-torsion group.

On the other hand, it is shown in~\cite{SaSt2}
that for an odd vector bundle $L$ on $S^{*}M$ and $N$ large enough
the bundle $2^NL$ is isomorphic to its complement $2^N\alpha ^{*}L$:
\[
\sigma _0:2^NL\stackrel\simeq\longrightarrow 2^N\alpha ^{*}L
\]
(this can be obtained noting that the
projection $S^{*}M\rightarrow P^{*}M$
for even-dimensional manifolds induces an isomorphism
in $K$-theory, modulo $2$-torsion). Now we can
construct the desired pull-back (\ref{omega}) by the formula
\[
\sigma :2^{N+1}L\stackrel{1\oplus \sigma _0}{\longrightarrow }2^NL\oplus
2^N\alpha ^{*}L=2^N\pi ^{*}E.
\]
\dokaend

Let us now consider an elliptic operator in subspaces
\[
\widehat{\sigma }:2^N\widehat{L}\longrightarrow C^\infty \left( M,F\right)
\]
with symbol (\ref{omega}).
We write out the index formula for this
operator:
\begin{equation}
\limfunc{ind}\left( \widehat{\sigma },2^N\widehat{L},C^\infty \left(
M,F\right) \right) =\frac 12\limfunc{ind}
\widehat{\widetilde{\sigma}}+2^Nd\left(
\widehat{L}\right) .  \label{omegasht}
\end{equation}
This formula, along with the integrality of the index, implies that the
functional $d$ is dyadic rational and has at most $2^{N+1}$ as the
denominator. For the fractional part of $d$ this gives
\[
\left\{d\left( \widehat{L}\right)\right\}=\frac 1{2^N}\left(
\mbox{mod $2^N$-ind}
\left( \widehat{\sigma },2^N\widehat{L},C^\infty \left( M,F\right)
\right) -\frac 12\left(\mbox{mod $2^{N+1}$-ind}
\widehat{\widetilde{\sigma }}\right)\right).
\]
In particular, for $N=0$ this gives a topological formula for the
fractional part of $d\left(\widehat{L}\right)$. For $N\ge 1$ it remains to
compute the index modulo $2^N$ for an elliptic operator in
subspaces
\[
\widehat{\sigma }:2^N\widehat{L}\longrightarrow C^\infty \left( M,F\right) .
\]
This problem is solved in the next section.

\section{Index theory modulo $\bf n$\label{comput}}

For a given positive integer $n\geq 2$, we consider elliptic operators in
subspaces of a special form:
\begin{equation}
D:n\widehat{L}\longrightarrow C^\infty \left( M,F\right) .  \label{lab3}
\end{equation}
Let us emphasize that the subspace $\widehat{L}$ need not satisfy the
parity condition. It is easy to show that the index of $D$ modulo $n$
denoted by
$$
\mbox{mod $n$-ind$\,D$} \in \Bbb{ Z}_n
$$
is determined by the principal symbol
$\sigma \left( D\right) :nL\longrightarrow \pi ^{*}F.$

It is natural to compute this index with values in
$\Bbb{ Z}_n,$ in terms of a difference construction with values in
$K$-theory with $\Bbb{ Z}_n$ coefficients:
\[
\left[\sigma \left( D\right) \right] \in K\left( T^{*}M,\Bbb{ Z}_n\right) .
\]
The necessary information about this theory is provided
for the reader's convenience in the Appendix.

Let us define this difference construction. First of all, we rewrite the
group $K\left( T^{*}M,\Bbb{ Z}_n\right) $ in terms of the usual $K$-theory.
We have
\begin{equation}
K\left( T^{*}M,\Bbb{ Z}_n\right) =K\left( T^{*}M\times \Bbb{M}_n,T^{*}M\times
pt\right) ,  \label{lab4}
\end{equation}
where $\Bbb{M}_n$ is the so-called {\em Moore space\/}. It readily follows from (\ref{lab4}) that the
elements of $K\left( T^{*}M,\Bbb{ Z}_n\right) $ can be realized
as families of elliptic symbols\footnote{Here we use the natural
construction \cite{AtSi4} that assigns an element
\[
\left[ \sigma \right] \in K\left( T^{*}M\times X\right)
\]
of the $K$-group to each
family $\sigma \left( x\right) ,$ $x\in X$, of elliptic symbols on the
manifold $M$ with the parameter space $X$:
\[
\sigma \left( x\right) :\pi ^{*}E\longrightarrow \pi ^{*}F,\quad E,F\in
\limfunc{Vect}\left( M\times X\right) ,\text{ }\pi :S^{*}M\times
X\rightarrow M\times X.
\]
} on the manifold $M,$ where the Moore space serves as the parameter space
for the family. Thus, the desired family of elliptic symbols is defined
as the composition of elliptic families in subspaces:
\begin{equation}
\begin{array}{c}
\left[ \sigma \left( D\right) \right] =
\left[ \pi ^{*}F\stackrel{\sigma ^{-1}\left( D\right) }{\longrightarrow }nL%
\stackrel{\beta ^{-1}\otimes 1_L}{\longrightarrow }\gamma _n\otimes nL%
\stackrel{1_\gamma \otimes \sigma \left( D\right) }{\longrightarrow }\gamma
_n\otimes \pi ^{*}F\right]
\vspace{3mm} \\
\left[ \sigma \left( D\right) \right] \in
K\left( T^{*}M\times \Bbb{M}_n,T^{*}M\times
pt\right) ,
\end{array}
\label{long}
\end{equation}
where $\gamma _n$ is the line bundle  corresponding to
the generator of the reduced group $\widetilde{K}\left( \Bbb{M}_n\right)
\simeq \Bbb{Z}_n $ and $\beta $ is a trivialization
$\beta :n\gamma _n\rightarrow \Bbb{ C}^n.$

\begin{theorem}
(index theorem modulo $n$)
\begin{equation}
\mbox{{\rm mod} $n$-{\rm ind}$\,D$}=p_{!}\left[ \sigma \left(
D\right) \right],   \label{lab5}
\end{equation}
where the direct image in $K$-theory (with coefficients)
\[
p_{!}:K\left( T^{*}M,\Bbb{ Z}_n\right) \longrightarrow K\left( pt,\Bbb{ Z}%
_n\right) =\Bbb{ Z}_n,
\]
is induced by the mapping $p:M\longrightarrow pt.$
\end{theorem}

\noindent {\em Proof. }Consider the following three families of elliptic
operators in subspaces, parametrized by the Moore space $\Bbb{M}_n$
\[
\begin{array}{c}
C^\infty \left( M,F\right) \stackrel{D^{-1}}{\longrightarrow }n\widehat{L}%
,\qquad \qquad \qquad \qquad \qquad \qquad \qquad  \\
\ n\widehat{L}\stackrel{\beta ^{-1}\otimes 1_{\widehat{L}}}{\longrightarrow }%
\gamma _n\otimes n\widehat{L}, \\
\qquad \qquad \qquad \qquad \qquad \qquad \qquad \gamma _n\otimes n\widehat{L%
}\stackrel{1_\gamma \otimes D}{\longrightarrow }\gamma _n\otimes C^\infty
\left( M,F\right)
\end{array}
\]
(here $D^{-1}$ denotes an almost inverse, i.e. inverse up to compact
operators, of $D$ and the three families correspond to the symbols in
(\ref{long})). The first family is constant. The second family is defined
by isomorphisms, while the third family is merely the
tensor product of $D$ and the bundle $\gamma _n$ over the
parameter space. Hence, the index of the composition  is equal to
\begin{equation}
\limfunc{ind}\left( \left[ 1_\gamma \otimes D\right] \circ \left[ \beta
^{-1}\otimes 1_{\widehat{L}}\right] \circ D^{-1}\right) =\left[ \gamma
_n\right] \limfunc{ind}D+0-\limfunc{ind}D\in K\left( \Bbb{M}_n\right) .
\label{lab6}
\end{equation}
On the other hand, the index of the elliptic family
\[
\left[ 1_\gamma \otimes D\right] \circ \left[ \beta ^{-1}\otimes 1_{\widehat{%
L}}\right] \circ D^{-1}:C^\infty \left( M,F\right) \longrightarrow \gamma
_n\otimes C^\infty \left( M,F\right)
\]
is calculated by the Atiyah-Singer index formula for families (see \cite
{AtSi4}). Thus, by virtue of (\ref{long}), this gives
\begin{equation}
\limfunc{ind}\left( \left[ 1_\gamma \otimes D\right] \circ \left[ \beta
^{-1}\otimes 1_{\widehat{L}}\right] \circ D^{-1}\right) =p_{!}\left[ \sigma
\left( D\right) \right] \in K\left( \Bbb{M}_n\right) .  \label{lab7}
\end{equation}
On the other hand, taking into account the isomorphism
$K\left( \Bbb{M}_n\right) =\Bbb{ Z\oplus Z}_n$ {with } $\left[ \gamma
_n\right] -1$ as a generator of the torsion part $\Bbb{ Z}_n,$
we obtain
\[
\mbox{mod $n$-ind$\,D$}=p_{!}\left[ \sigma \left(
D\right) \right]
\]
comparing (\ref{lab7}) with (\ref{lab6}).
\dokaend

\begin{remark}
{\em A similar mod$n$-index theorem for boundary value problems
was obtained in \cite{Fre1}, \cite{FrMe1}.}
\end{remark}

\section{A formula for the fractional part}

In this section we write out the final formula for the $\eta$-invariant.
Namely, for a subspace
$\widehat{L}\subset C^\infty \left( M,E\right) $ and the pull-back of its
symbol from $M$ by an isomorphism $\sigma $
$$
2^NL\stackrel{\sigma }{\longrightarrow }\pi ^{*}F,\quad \quad F\in
\limfunc{Vect}\left( M\right) ,\pi :S^{*}M\rightarrow M,
$$
in Section \ref{secsimp} we expressed the fractional part
of $d(\widehat{L})$ as
\[
\left\{d\left( \widehat{L}\right)\right\}=\frac 1{2^N}\left(
\mbox{mod $2^N$-ind}
\left( \widehat{\sigma },2^N\widehat{L},C^\infty \left( M,F\right)
\right) -\frac 12\left(\mbox{mod $2^{N+1}$-ind}
\widehat{\widetilde{\sigma }}\right)\right).
\]
The two terms, in fact, can be combined together. Namely, in
the even (odd) cases the resulting formulas are, respectively,
\begin{equation}
\begin{array}{l}
\left\{d\left( \widehat{L}\right)\right\}= \\
\frac 1{2^{N+1}}\mbox{mod $2^{N+1}$-ind}
\left( 2^{N+1}\widehat{L}\stackrel{\widehat{%
\sigma }\oplus \widehat{\alpha ^{*}\sigma }}{\longrightarrow }C^\infty
\left( M,F\oplus F\right) \right)
,\quad \widehat{L}\in
\widehat{\limfunc{Even}}\left( M^{odd}\right) , \\
\frac 1{2^{N+1}}\mbox{mod $2^{N+1}$-ind}
\left( 2^{N+1}\widehat{L}\stackrel{1\oplus
\widehat{\alpha ^{*}\left[ \sigma ^{-1}\right] }\widehat{\sigma }}{%
\longrightarrow }C^\infty \left( M,2^NE\right) \right),\quad \widehat{L}\in \widehat{\limfunc{Odd}}\left( M^{ev}\right) .
\end{array}
\label{evad}
\end{equation}
Applying the mod$n$-index formula, we obtain the desired
topological expression
\begin{equation}
\label{infk}
2^{N+1}\{ d( \widehat{L})\}
=p_{!}\left[ L\right]_N \in \Bbb{ Z}_{2^{N+1}},\quad \left[ L\right]_N \in
K\left( T^{*}M,\Bbb{ Z}_{2^{N+1}}\right) ,
\end{equation}
where $\left[ L\right]_N $ denotes the
difference construction for the operators in (\ref{evad})
\begin{equation}
\label{fina}
\left[ L\right] =\left[ 2^{N+1}L\stackrel{\sigma \oplus \alpha ^{*}\sigma }{%
\longrightarrow }\pi ^{*}F\oplus \pi ^{*}F\right] \text{\quad or\quad }%
\left[ 2^{N+1}L\stackrel{1\oplus \alpha ^{*}\left[ \sigma ^{-1}\right]
\sigma }{\longrightarrow }2^N\pi ^{*}E\right] .
\end{equation}
We would like to rewrite (\ref{infk}) in a more canonical form.

To this end, consider the embedding
$$
i:\Bbb{Z}_{2^{N+1}}\subset\Bbb{Z}\left.\left[\frac 12\right]
\right/\Bbb{Z}, \quad i(x)=
\frac{x}{2^{N+1}}.
$$
It induces a mapping of $K$-groups
$$
i_*:K\left(T^*M,\Bbb{Z}_{2^{N+1}}\right)\longrightarrow
K\left(T^*M,\left.\Bbb{Z}\left[\frac 12\right]\right/\Bbb{Z}\right),
$$
where the $K$-theory with dyadic coefficients is defined as the injective
limit
\begin{equation}
K\left( T^{*}M,\left. \Bbb{ Z}\left[ \frac 12\right] \right/ \Bbb{ Z}\right)
=\lim_{
\scriptstyle \overrightarrow{
\scriptstyle N^{\prime }\rightarrow \infty }}K\left( T^{*}M,\Bbb{ Z}%
_{2^{N^{^{\prime }}}}\right) .  \label{lima}
\end{equation}

\begin{lemma}
\label{qq}The element $\left[ L\right]=i_*[L]_N \in
K\left( T^{*}M,\left. \Bbb{ Z}\left[1/2\right] \right/ \Bbb{ Z}\right)$
is well defined, i.e. independent of the choice of isomorphism $\sigma $.
\end{lemma}

\noindent {\em Proof.} For two isomorphisms
\[
2^NL\stackrel{\sigma }{\longrightarrow }\pi ^{*}F,\quad \text{and\quad }2^NL%
\stackrel{\sigma ^{\prime }}{\longrightarrow }\pi ^{*}F^{\prime }
\]
let us compute the difference of the corresponding $K$-theory elements in (%
\ref{fina}). An explicit computation shows that the difference is equal to
\begin{eqnarray*}
&&\left[\sigma \sigma ^{\prime -1}\right]
\oplus \alpha ^{*}\left[ \sigma \sigma ^{\prime
-1}\right] \qquad \text{on an odd-dimensional manifold }M, \\
&&
\left[\sigma \sigma ^{\prime -1}\right]\oplus \alpha ^{*}\left[ \sigma ^{\prime }\sigma
^{-1}\right] \text{\qquad on an even-dimensional }M.
\end{eqnarray*}

Thus, the difference in question is equal to
\[
\left[ \sigma _0\right] \pm \left[ \alpha ^{*}\sigma _0\right] \in K\left(
T^{*}M\right)
\]
for the elliptic symbol $\sigma _0=\sigma \sigma ^{\prime -1}.$
The sign is opposite to the parity of ${\rm dim}M.$

To prove the Lemma, it suffices to show that $\sigma_0$
defines a 2-torsion element in  $K\left(T^{*}M\right).$
This is proved in the following purely topological Proposition.
\dokaend

\begin{proposition}
\label{th5}
The antipodal involution $\alpha :T^{*}M\longrightarrow T^{*}M$
induces an involution
\[
\alpha ^{*}:K^*\left( T^{*}M\right) \longrightarrow K^*\left( T^{*}M\right)
\]
equal to $\left( -1\right) ^{\dim M}$, modulo $2$-torsion. More
precisely, for an arbitrary $x\in K\left( T^{*}M\right) $ and
$N$ large enough one has
\begin{equation}
\alpha ^{*}\left( 2^Nx\right) =\left( -1\right) ^{\dim M}2^Nx.  \label{power}
\end{equation}
\end{proposition}

\noindent {\em Proof. }The idea is to apply the
Mayer--Vietoris principle.

1) Let us prove ($\ref{power}$) for $M=\Bbb{R}^n$. We have
\[
K^{*}\left( T^{*}\Bbb{R}^n\right) \simeq K^{*}\left( \Bbb{ R}^{2n}\right)
=\Bbb{Z},
\]
while $\alpha:\Bbb{R}^{2n}\to \Bbb{R}^{2n}$ acts as
$(x,\xi)\to(x,-\xi)$. Thus, for $n$ even it is homotopic
to the identity and in the $K$-theory we have $\alpha^*=id$.
While in odd-dimensions $\alpha$ reverses the orientation.
Therefore, in this case $\alpha^*=-id$, as desired.

2) We claim that the following assertion is valid:
suppose that (\ref{power}) is satisfied for two open subsets
$U,V\subset M$ and for their intersection $U\cap V$. Then these properties
are valid for the union $U\cup V$.%

A part of the Mayer-Vietoris exact sequence
corresponding to the inclusions $U\bigcap V\stackrel{i}{\subset }U\sqcup V%
\stackrel{j}{\subset }U\bigcup V$ is
\[
\begin{array}{ccccc}
K^{*+1}\left( T^{*}\left( U\bigcap V\right) \right) & \stackrel{\delta }{%
\rightarrow } & K^{*}\left( T^{*}\left( U\bigcup V\right) \right) &
\stackrel{j^{*}}{\rightarrow } & K^{*}\left( T^{*}U\right) \oplus
K^{*}\left( T^{*}V\right) \\
\downarrow \alpha ^{*} &  & \downarrow \alpha ^{*} &  & \quad \quad \quad
\downarrow \!\!\alpha ^{*}\!\!\oplus \!\alpha ^{*} \\
K^{*+1}\left( T^{*}\left( U\bigcap V\right) \right) & \stackrel{\delta }{%
\rightarrow } & K^{*}\left( T^{*}\left( U\bigcup V\right) \right) &
\stackrel{j^{*}}{\rightarrow } & K^{*}\left( T^{*}U\right) \oplus
K^{*}\left( T^{*}V\right) .
\end{array}
\]
Suppose that the left and the right involutions in the diagram are equal to
$\left( -1\right) ^{\dim M}$ (modulo 2-torsion).
By a diagram chasing argument one deduces that the mapping
$\alpha ^{*}$ in the center also satisfies (\ref{power}).
For example, on an even-dimensional
manifold for $x\in K^{*}\left( T^{*}\left(
U\bigcup V\right) \right) $ we get
\[
j^{*}\left( \alpha ^{*}x-x\right) =0\Rightarrow \alpha ^{*}x-x=\delta \alpha
^{*}y,\alpha ^{*}y=y\Rightarrow 2\left( \alpha ^{*}x-x\right) =0
\]
(in this computation factors $2^N$ are omitted for brevity).

3) Consider a good (see \cite{BoTu1}) finite covering
$\left\{ U_\beta \right\} $
of the manifold $M$ by contractible open sets. Over any
$U_\beta $ the property (\ref{power}) is valid by the first part of the proof.
Let us consider all subsets in $\left\{ U_\beta \right\} .$

Passing from the coverings consisting of a single element to the covering
of the entire manifold $M$ and applying the assertion from the second part
of the proof, we obtain the desired property for $M$.

Now Eq. (\ref{infk}) and Lemma \ref{qq} prove the following theorem.
\begin{theorem}
\label{lab10}
A subspace $\widehat{L}\in \widehat{\limfunc{Even}}%
\left( M^{odd}\right) $ or $\widehat{\limfunc{Odd}}\left( M^{ev}\right) $
defines an element
\[
\left[ L\right] \in K\left( T^{*}M,\left. \Bbb{ Z}\left[ \frac 12\right]
\right/ \Bbb{ Z}\right),
\]
and the fractional part of the functional $d$ is computed by the direct
image mapping
$$
\left\{d\left( \widehat{L}\right)\right\}=p_{!}\left[ L\right] \in
K\left( pt,\left. \Bbb{ Z}\left[ \frac 12\right] \right/ \Bbb{ Z}\right)
=\left. \Bbb{ Z}\left[ \frac 12\right] \right/ \Bbb{ Z}
$$
induced by $p:\,M\rightarrow pt$. In terms of an isomorphism
$\sigma :2^NL\longrightarrow \pi ^{*}F$
(see Proposition {\em \ref{thpow}}), $\left[ L \right] $ is defined
by the symbol $\left[ \left( 1\pm \alpha
^{*}\right) \sigma \right] \in K\left( T^{*}M,\Bbb{ Z}_{2^{N+1}}\right) $
as
\[
\left[ L \right] =i_{*}\left[ \left( 1\pm \alpha
^{*}\right) \sigma \right] ,\quad i:\Bbb{ Z}_{2^{N+1}}\subset \Bbb{ Z}%
\left[ \frac 12\right] /\Bbb{ Z}
\]
(the sign coincides with the parity of the subspace).
A similar formula holds for the $\eta$-invariant
$$
\{\eta(A)\}=p_![L_+(A)]
$$
of an elliptic self-adjoint differential operator $A$ satisfying the
parity condition {\em (\ref{conda})}.
\end{theorem}

\section{Examples and remarks}
{\bf 1. Operators from} \cite{Gil7}.
On a smooth oriented closed Riemannian odd-dimensional manifold $M$, we
consider an elliptic self-adjoint differential operator of the second order
\begin{equation}
\label{aloe}
A=d\delta -\delta d:C^\infty \left( M,\Lambda ^1\left( M\right) \right)
\longrightarrow C^\infty \left( M,\Lambda ^1\left( M\right) \right)
\end{equation}
in the spaces of complexified 1-forms; here $d$ is the exterior derivative
and $\delta$ is the adjoint operator with respect to the Riemannian metric.
The principal symbol of $A$ is
\[
\sigma \left( A\right) \left( \xi \right) =\xi \wedge \xi \rfloor -\xi
\rfloor \xi \wedge :\pi ^{*}\Lambda ^1\left( M\right) \longrightarrow \pi
^{*}\Lambda ^1\left( M\right) ,
\]
where $\xi \wedge $ is the exterior product by a covector
$\xi$ and $\xi \rfloor $ is the inner product by the same covector with
respect to the Riemannian metric.
For an arbitrary point $\left( x,\xi \right) \in S_x^{*}M$
of the cosphere bundle, the symbol $L=L_+(A)$ of the spectral subspace
$\widehat L=\widehat{L}_+(A)$ coincides with the line spanned by the covector
$\xi $. Hence, $L\subset \pi ^{*}\Lambda ^1\left( M\right) $ is an
even subbundle that at $x\in M$ generates the reduced $K$-group of the
projective space
$$
[L]-[1]\in\widetilde{K}(P_x^*M)\simeq\Bbb{Z}_{2^{({\rm dim}M-1)/2}},
\quad \mbox{ for}\;\;{\rm dim}M\ge 5.
$$
Thus, the operator $A$ defines a nontrivial element of $K(P^*M)/K(M)$.
Let us compute the fractional part of the functional $d$ on the
subspace $\widehat{L}_+(A)$.

The line bundle $L$ is trivial. We choose the trivialization
\begin{equation}
\begin{array}{ccc}
\sigma :L & \rightarrow  & \pi ^{*}\Bbb{C}, \\
\sigma \left( x,\xi \right) \eta  & = & \left\langle \xi ,
\eta \right\rangle_x
\end{array}
,\qquad \left( x,\xi \right) \in S_x^{*}M,\;\eta \in L_x,  \label{simple}
\end{equation}
where $\left\langle \xi ,\eta \right\rangle_x $ denotes the Hermitian
inner
product of two covectors with respect to the Riemannian metric at  $x$.
For the corresponding pseudodifferential operator
\[
\widehat{\sigma }:\widehat{L}\longrightarrow C^\infty \left( M\right)
\]
in the subspaces, the index formula (\ref{inde}) implies
\[
\limfunc{ind}\left( \widehat{\sigma },\widehat{L},C^\infty \left( M\right)
\right) =\frac 12\limfunc{ind}\widehat{\widetilde{\sigma }}+d\left( \widehat{%
L}\right) .
\]
It follows from (\ref{simple}) that the symbol $\widetilde{%
\sigma }$ is a direct sum of two constant symbols
\[
\begin{array}{c}
\widetilde{\sigma }:\pi ^{*}\Lambda ^1\left( M\right) \rightarrow \pi
^{*}\Lambda ^1\left( M\right) , \\
\widetilde{\sigma }\left( \xi \right) =\sigma ^{-1}\left( -\xi \right)
\sigma \left( \xi \right) \oplus 1=-1\oplus 1:L\oplus L^{\perp }\rightarrow
L\oplus L^{\perp }.
\end{array}
\]
Hence, we obtain the integrality result for the functional $d$:
\[
d\left( \widehat{L}\right) =\limfunc{ind}\left( \widehat{\sigma },\widehat{L}%
,C^\infty \left( M\right) \right) \in \Bbb{Z}.
\]

By the same method one can show the integrality of the functional for
operators $A$
with coefficients in a vector{\rm \ }bundle $E\in \limfunc{Vect}\left(
M\right) $. To this end, one replaces the exterior derivative
$d$ and its adjoint by a covariant derivative $\nabla $ and the
corresponding adjoint operator for $E$.
These operators were introduced in \cite{Gil7}, where the problem of
nontriviality of their $\eta$-invariants was posed. Thus, we obtain
\begin{proposition}
\label{Gi1}
Operator $d\delta-\delta d$ with coefficients in a bundle $E$ has
trivial fractional part of the $\eta$-invariant.
\end{proposition}
As an application, consider the $3$-dimensional real projective space
$\Bbb{RP}^3$.  It is parallelizable and the Kunneth formula shows that
the group $K(P^*\Bbb{RP}^3)/K(\Bbb{RP}^3)$ is generated by operators
(\ref{aloe}) with coefficients in vector bundles. Thus, this manifold
has no even order operators with fractional $\eta$-invariants.

{\bf 2. } The vanishing result of Proposition \ref{Gi1} has an interesting
corollary.

Assume $M$ is odd-dimensional as before. Consider a pair $S^*M\subset B^*M$,
where both spaces have the natural
antipodal $\Bbb{Z}_2$ action. The equivariant exact sequence of this pair
leads to the following long exact sequence
\begin{equation}
\!\to\! K(M)\!\rightarrow K(P^*M)/K(M)\rightarrow K^1_{\Bbb{Z}_2}(T^*M)\!
\to K^1(M)\! \to \!K^1(P^*M)/K^1(M)\!\to\ldots
\end{equation}
In analytic terms, the mapping $K(M)\to K(P^*M)/K(M)$ corresponds to
taking a bundle $E$ to (symbol of) the operator $d\delta-\delta d$ with
coefficients in $E$.

Thus, the $\eta$-invariant as a mapping
$K(P^*M)/K(M)\to \Bbb{Z}[1/2]/\Bbb{Z}$, is induced by a mapping
of the subgroup in $K^1_{\Bbb{Z}_2}(T^*M)$.
The elements of this latter equivariant $K$-group can be realized as symbols
of indicial families in the sense of \cite{Mel6} with a special symmetry.
It seems that the $\eta$-invariant of Melrose would give the
corresponding analytic realization of this mapping.
However, this topic is beyond the scope of the present paper.

\newpage

{\bf 3. Orientable manifolds and the $\eta$-invariant.}
\begin{proposition}
If $M$ is orientable then the functional $d$ and the $\eta$-invariant
take integer, or half-integer values only.
\end{proposition}

\noindent
{\em Proof.} Indeed, Theorem \ref{lab10} gives
$$
2^{N+1}\{d(\widehat{L})\}=p_![(1-(-1)^{{\rm dim}M}\alpha^*)\sigma]\in
\Bbb{Z}_{2^{N+1}}.
$$
Thus,
$$
2^{N+2}\{d(\widehat{L})\}=2p_!(1-(-1)^{{\rm dim}M}\alpha^*)
[\sigma\oplus\sigma]\in
\Bbb{Z}_{2^{N+1}},\quad
[\sigma\oplus\sigma]\in K(T^*M,\Bbb{Z}_{2^{N+1}}).
$$
However, M.~Karoubi proved in \cite{Karo3} that on an orientable
$M$ the antipodal involution $\alpha$ acts as
 $(-1)^{{\rm dim}M}$ in $K$-theory. Therefore, this yields
$$
2^{N+2}\{d(\widehat{L})\}=p_!0=0.
$$
Hence, we prove the desired ${2d(\widehat{L})}\in \Bbb{Z}.$
\dokaend

\section{Elliptic theory with $\Bbb{Z}_{n}$ coefficients\label{elln}}

The difference construction (\ref{long}) is not an entirely computational
trick involved in the modulo $n$-index calculation above. In the present
section we show that similar to the usual difference construction,
it establishes an isomorphism between the $K$-theory
$K\left( T^{*}M,\Bbb{ Z}_n\right) $ with $\Bbb{ Z}_n$ coefficients and the
group of stable homotopy classes of elliptic operators in subspaces of the
form (\ref{lab3}).

{\bf 1. Definition.} We consider elliptic operators of the form
\begin{equation}
D=n\widehat{L}_1\oplus C^\infty \left( M,E_1\right) \longrightarrow n%
\widehat{L}_2\oplus C^\infty \left( M,F_1\right) ,
\label{lab8}
\end{equation}
where $\widehat{L}_1\subset C^\infty \left( M,E\right),
\widehat{L}_2\subset C^\infty \left( M,F\right)$.
This is slightly different from (\ref{lab3}); the difference is
motivated by the requirement that the inverse operator be of the same
structure.

Let us state the stable homotopy classification problem for
such operators. First, we introduce {\em trivial} operators. These
are: a) operators induced by  a vector bundle isomorphisms
$g:E_1\rightarrow F_1$:
\begin{equation}
C^\infty \left( M,E_1\right) \stackrel{g_{*}}{\longrightarrow }C^\infty
\left( M,F_1\right),   \label{one}
\end{equation}
b) direct sums of $n$ copies of an elliptic operator in subspaces:
\begin{equation}
n\left( \widehat{L}_1\oplus C^\infty \left( M,E_1\right) \right) \stackrel{nD%
}{\longrightarrow }n\left( \widehat{L}_2\oplus C^\infty \left( M,F_1\right)
\right) .  \label{two}
\end{equation}

We identify operators of the form (\ref{lab8}) that differ by isomorphisms
of the corresponding vector bundles $E,F,E_1,F_1.$ Two elliptic operators $%
D_1$ and $D_2$ are {\em stably homotopic} if they become homotopic after
we add some trivial operators to each of them. The
abelian group formed by
the classes of stably homotopic elliptic operators
is denoted by $\limfunc{Ell}\left( M,\Bbb{ Z}_n\right) .$
As usual, one can prove that the composition of elliptic operators $D_{1,2}$
(if defined) gives an element $[D_1D_2]$
equal to the sum $[D_1]+[D_2]$.

\begin{lemma}
\label{lemm}An operator {\em (\ref{lab8})} is stably homotopic
to an operator of the form
\begin{equation}
n\widehat{L}'\stackrel{D'}{\longrightarrow }C^\infty \left(
M,F'\right).
\label{simp1}
\end{equation}
\end{lemma}

\noindent {\em Proof. }The space $C^\infty \left( M,E_1\right) $ can be
eliminated in (\ref{lab8}) by adding the trivial operator
\[
\left( n-1\right) \left( C^\infty \left( M,E_1\right) \stackrel{id}{%
\longrightarrow }C^\infty \left( M,E_1\right) \right) .
\]
The subspace $\widehat{L}_2$ on the right-hand side of the formula can be
eliminated in the following way. Let us add the trivial operator
$id:n\widehat{L}_2^{\perp }\to n\widehat{L}_2^{\perp }$ to $D$.
Then we obtain an operator of the form
\[
n\left( \widehat{L}_1\oplus \widehat{L}_2^{\perp }\right) \longrightarrow
n\left( \widehat{L}_2\oplus \widehat{L}_2^{\perp }\right) \oplus C^\infty
\left( M,F_1\right).
\]
To complete the proof, it suffices to show that the
resulting subspace
\[
\widehat{L}_2\oplus \widehat{L}_2^{\perp }\subset C^\infty \left(
M,F_1\right) \oplus C^\infty \left( M,F_1\right)
\]
is homotopic to the subspace $C^\infty \left( M,F_1\right) \oplus 0$,
since a homotopy of subspaces can be lifted to a homotopy of elliptic
operators. The desired homotopy of subspaces is given in terms
of the projection $P$ on the subspace $\widehat{L}_2$ by the
formula
\[
\widehat{L}_\varphi =\func{Im}P_\varphi ,\;P_\varphi =\left(
\begin{array}{cc}
P & 0 \\
0 & 0
\end{array}
\right) +\left( 1-P\right) \left(
\begin{array}{cc}
\sin ^2\varphi & \cos \varphi \sin \varphi \\
\cos \varphi \sin \varphi & \cos ^2\varphi
\end{array}
\right) .
\]
Here $\widehat{L}_\varphi \subset C^\infty \left( M,F_1\right) \oplus
C^\infty \left( M,F_1\right) $ and $\varphi $ varies
from $0$ to $\pi /2.$
\dokaend

{\bf 2. Exact sequence in Elliptic theory.}
Denote by $\limfunc{Ell}\left( M\right) $
the group of stable homotopy classes of elliptic operators on $M$.
Let $\limfunc{Ell}_1\left( M\right) $ denote a similar group
of stable homotopy classes of pseudodifferential subspaces.
More precisely, two subspaces are called homotopic, if
there is a norm continuous homotopy of projections defining them.
They are stably homotopic, if they become homotopic if
we add some {\em trivial subspaces} to them. Here the trivial subspaces
are spaces of sections of vector bundles and finite-dimensional spaces.

Let us now define a sequence
\begin{equation}
\limfunc{Ell}\left( M\right) \stackrel{\times n}{\longrightarrow }\limfunc{%
Ell}\left( M\right) \stackrel{i}{\longrightarrow }\limfunc{Ell}\left( M,\Bbb{ %
Z}_n\right) \stackrel{j}{\longrightarrow }\limfunc{Ell}_1\left(
M\right) \stackrel{\times n}{\longrightarrow }\limfunc{Ell}_1\left(
M\right) ,  \label{longtwo}
\end{equation}
where $\times n$ denotes multiplication by $n$, the mapping $i$ is induced
by the natural inclusion of the usual elliptic operators in
the mod$n$-theory, and the boundary mapping $j$ is defined as
\[
j\left[ n\widehat{L}_1\oplus C^\infty \left( M,E_1\right) \stackrel{D}{%
\longrightarrow }n\widehat{L}_2\oplus C^\infty \left( M,F_1\right) \right]
=\left[ \widehat{L}_1\right] -\left[ \widehat{L}_2\right] .
\]

\begin{proposition}
{\em The sequence (\ref{longtwo}) is exact.}
\end{proposition}

\noindent {\em Proof. }It is straightforward to show that (\ref{longtwo})
is a complex. Let us prove the exactness.

Let $\left[ D\right] \in \ker j.$ By Lemma \ref{lemm} we can suppose that
$D$ has the form (\ref{simp1}). Since $\left[ \widehat{L}'
\right] =0\in \limfunc{Ell}_1\left( M\right) ,$ it follows that
the subspace $\widehat{L}'$ is homotopic\footnote{%
Here and in what follows we omit the standard considerations concerning the
stabilization of elements.} to a space of sections
of a vector bundle. Hence, $D$ is homotopic
to an elliptic operator in the spaces of vector bundle sections.
Consequently, we obtain $\left[ D\right] \in \func{Im}i.$

Let $\left[ \widehat{L}\right] \in \ker \left\{ \times n\right\}.$
This implies that the subspace $n%
\widehat{L}$ is homotopic to the space of sections of a vector bundle.
Consequently, there exists an elliptic operator
$D:n\widehat{L}{\rightarrow }C^\infty \left( M,F\right) .$
Hence, $[ \widehat{L}] =j\left[ D\right],$ as desired.
The remaining assertion $\left[ D\right] \in \ker i$ $\Rightarrow$
$[D]\in{\rm Im}(\times n)$
can be proved along similar lines and is left to the reader.
\dokaend

{\bf 3. Isomorphism of Elliptic theory and $K$-theory.}
By virtue of Lemma \ref{simp1}, we can extend the difference construction
\begin{equation}
D\longmapsto \left[ \sigma \left( D\right) \right] \in K\left( T^{*}M,\Bbb{ Z}%
_n\right)  \label{twostar}
\end{equation}
(see (\ref{long})) to a homomorphism of groups
\[
\limfunc{Ell}\left( M,\Bbb{ Z}_n\right) \longrightarrow K\left( T^{*}M,\Bbb{ Z}%
_n\right) ,
\]
since the mapping (\ref{long}) sends the trivial operators (\ref{one}) and (\ref
{two}) to zero.

Let us also recall the difference construction for pseudodifferential
subspaces. Namely, a subspace $\widehat{L}=\func{Im%
}P$ with symbol $L=\func{Im}\sigma \left( P\right) $ defines
a family of elliptic symbols on $M$ with the
parameter space $\Bbb{S}^1$ and coordinate $z$:
\[
z\sigma \left( P\right) +\left( 1-\sigma \left( P\right) \right) :\pi
^{*}E\longrightarrow \pi ^{*}E.
\]
By virtue of the usual difference construction for elliptic families,
this defines the desired element (see \cite{APS3})
\begin{equation}
\label{difodd}
\left[ z\sigma \left( P\right) +\left( 1-\sigma \left( P\right) \right)
\right] \in K\left( T^{*}M\times \Bbb{S}^1,T^{*}M\times pt\right) \equiv K^1\left(
T^{*}M\right).
\end{equation}

Let us now consider the following diagram
\begin{equation}
\begin{array}{ccccccccc}
\limfunc{Ell}\left( M\right) & \stackrel{\times n}{\!\!\rightarrow\!\! } & \limfunc{%
Ell}\left( M\right) & \stackrel{i}{\!\!\rightarrow\!\! } & \limfunc{Ell}\left( M,%
\Bbb{ Z}_n\right) & \stackrel{j}{\!\!\rightarrow\!\! } & \limfunc{Ell}%
_1\left( M\right) & \stackrel{\times n}{\!\!\rightarrow\!\! } & \limfunc{Ell%
}_1\left( M\right) \\
\downarrow \chi _0 &  & \downarrow \chi _0 &  & \downarrow \chi _n &  &
\downarrow \chi _1 &  & \downarrow \chi _1 \\
K\left( T^{*}M\right) & \stackrel{\times n}{\!\!\rightarrow\!\! } & K\left(
T^{*}M\right) & \stackrel{i^{\prime }}{\!\!\rightarrow\!\! } & K\left( T^{*}M,\Bbb{ Z}%
_n\right) & \stackrel{j^{\prime }}{\!\!\rightarrow\!\! } & K^1\left( T^{*}M\right) &
\stackrel{\times n}{\!\!\rightarrow\!\! } & K^1\left( T^{*}M\right),
\end{array}
\label{vlong}
\end{equation}
where $\chi $ with subscripts denote difference constructions, and the
lower row in the diagram is part of the exact coefficient sequence in
$K$-theory (see Appendix).

\begin{lemma}
\label{seben}The diagram{\em \ (\ref{vlong}) } is commutative.
\end{lemma}

\noindent {\em Proof . }The commutativity of the
leftmost and rightmost squares of the diagram is clear. Let us
consider the second square:
\[
\begin{array}{ccc}
\limfunc{Ell}\left( M\right) & \stackrel{i}{\longrightarrow } & \limfunc{Ell}%
\left( M,\Bbb{ Z}_n\right) \\
\downarrow \chi _0 &  & \downarrow \chi _n \\
K\left( T^{*}M\right) & \stackrel{i^{\prime }}{\longrightarrow } & K\left(
T^{*}M,\Bbb{ Z}_n\right) .
\end{array}
\]
For an elliptic operator $D$ it is easy to see that
\[
\chi _ni\left[ D\right] =\left[ \sigma \left( D\right) \right] \left( \left[
\gamma _n\right] -1\right) \in K\left( T^{*}M\times \Bbb{M}_n,T^{*}M\times
pt\right) ;
\]
here $\left[ \sigma \left( D\right) \right] =\chi _0\left[ D\right] \in
K\left( T^{*}M\right) $ is the usual difference construction. On the other
hand, the reduction modulo $n$ mapping $i^{\prime }$ is exactly the
multiplication by the element $\left[ \gamma _n\right] -1.$ Thus, the second
square in (\ref{vlong}) is commutative.

Finally, let us check the commutativity of the remaining third square
\[
\begin{array}{ccc}
\limfunc{Ell}\left( M,\Bbb{ Z}_n\right) & \stackrel{j}{\longrightarrow } &
\limfunc{Ell}_1\left( M\right) \\
\downarrow \chi _n &  & \downarrow \chi _1 \\
K\left( T^{*}M,\Bbb{ Z}_n\right) & \stackrel{j^{\prime }}{\longrightarrow } &
K^1\left( T^{*}M\right) .
\end{array}
\]
For an elliptic operator
$D:n\widehat{L}{\longrightarrow }C^\infty \left( M,F\right),$
on the one hand, we obtain
\[
\chi _1\left( j\left[ D\right] \right) =\left[ L\right] \in K^1\left(
T^{*}M\right) .
\]
On the other hand, the difference construction for $D$ gives
\begin{equation}
\chi _n\left[ D\right] =\left[ \pi ^{*}F\stackrel{\sigma ^{-1}\left(
D\right) }{\longrightarrow }nL\stackrel{\beta ^{-1}\otimes 1_L}{%
\longrightarrow }\gamma _n\otimes nL\stackrel{1_\gamma \otimes \sigma \left(
D\right) }{\longrightarrow }\gamma _n\otimes \pi ^{*}F\right] \in K\left(
T^{*}M,\Bbb{ Z}_n\right) .  \label{llong}
\end{equation}
In terms of the identifications
\begin{eqnarray*}
K\left( T^{*}M,\Bbb{ Z}_n\right) &=&K\left( T^{*}M\times \Bbb{M}_n,T^{*}M\times
pt\right) , \\
K^1\left( T^{*}M\right) &=&K\left( T^{*}M\times \Bbb{S}^1,T^{*}M\times pt\right) ,
\end{eqnarray*}
the Bokstein mapping $j^{\prime }$
is induced by the inclusion $\Bbb{S}^1\stackrel{i_0}{\subset }\Bbb{M}_n$. More
precisely,
\[
j^{\prime }=\left( 1_{T^{*}M}\times i_0\right) ^{*}.
\]

Let us compute the family of elliptic symbols in (\ref{llong}) on the
circle $\Bbb{S}^1\subset \Bbb{M}_n$ with a polar coordinate
$\zeta =e^{i\varphi },0\leq \varphi <2\pi$. The family in (\ref{llong})
has the form
\begin{equation}
\pi ^{*}F\stackrel{\sigma ^{-1}}{\longrightarrow }nL\stackrel{\zeta
\oplus 1_{n-1}
}{\longrightarrow }nL\stackrel{\sigma }{\longrightarrow }\pi ^{*}F
\label{shtsht}
\end{equation}
with respect to the natural trivialization of $\gamma _n$ on $\Bbb{S}^1$.
Here the principal symbol of $D$ is denoted by $\sigma $ and
the diagonal operator $\zeta \oplus 1$ acts as
$\left( \zeta \oplus 1\right) \left( u_1,u_2,\ldots ,u_n\right)
=\left( \zeta u_1,u_2,\ldots ,u_n\right) .$

Let us also rewrite the symbol $\sigma $ in block matrix form:
\[
\sigma =\left( \sigma _1,\ldots ,\sigma _n\right) ,\qquad \sigma
_i:L\longrightarrow \pi ^{*}F.
\]
The ellipticity of $\sigma $ implies that the components $%
\sigma _i$ are monomorphic. Consider also the inverse symbol $\sigma ^{-1}$%
\[
\sigma ^{-1}=\left(\sigma ^1 ,\cdots,\sigma ^n\right)^t ,
\qquad \sigma ^i:\pi ^{*}F\longrightarrow L.
\]
We readily obtain
\[
\sum_{i=1}^n\sigma _i\sigma ^i=1,\qquad \sigma ^i\sigma _j=\delta _j^i.
\]
Hence, $\sigma _1\sigma ^1$ is the projection on a subbundle
isomorphic to the original bundle $L$%
\[
\func{Im}\sigma _1\sigma ^1\stackrel{\sigma _1}{\simeq }L.
\]
Therefore, the family (\ref{shtsht}) defines an element
\[
[\zeta \sigma _1\sigma ^1+\left( 1-\sigma _1\sigma ^1\right)] \in K\left(
T^{*}M\times \Bbb{S}^1,T^{*}M\times pt\right) =K^1\left( T^{*}M\right).
\]
This element coincides with the difference construction for
$\widehat{L}$ (see (\ref{difodd})).
\dokaend

\begin{theorem}
The difference construction
\[
\limfunc{Ell}\left( M,\Bbb{ Z}_n\right) \stackrel{\chi _n}{\longrightarrow }%
K\left( T^{*}M,\Bbb{ Z}_n\right)
\]
is an isomorphism.
\end{theorem}

\noindent
{\em Proof.} The usual difference constructions
\[
\limfunc{Ell}\left( M\right) \stackrel{\chi _0}{\rightarrow }K\left(
T^{*}M\right) \quad \text{and \quad }\limfunc{Ell}_1\left( M\right)
\stackrel{\chi _1}{\rightarrow }K^1\left( T^{*}M\right)
\]
are isomorphisms, so the theorem is proved by applying the 5-lemma to the
commutative diagram (\ref{vlong}).
\dokaend

\appendix
\vspace{0.5cm}

\noindent
{\LARGE\bf Appendix. $K$-theory with
coefficients}
\addcontentsline{toc}{section}{Appendix. $K$-theory with coefficients}

Here we recall some basic properties of the $K$-theory with
$\Bbb{ Z}_n$ coefficients that are used in the present paper.
More details can be found, e.g. in \cite{ArTo1,Bla1},
and the references therein.
By $n$ we denote a positive integer, $n\ge 2$.

{\bf 1. Moore space.} Let us consider the $2$-dimensional complex $
\Bbb{M}_n $ obtained from the unit disk $D^2$ identifying points on its
boundary under the $\Bbb{ Z}_n$ action:
\[
\Bbb{M}_n=\left. \left\{ \left. D^2\subset \Bbb{ C}\right| \;\left| z\right| \leq
1\right\} \right/ \left\{ e^{i\varphi }\sim e^{i\left( \varphi +\frac{2\pi k}%
n\right) }\right\} .
\]
The result is called the {\em Moore space. }For instance, $\Bbb{M}_2=\Bbb{ RP}^2$.
There is an embedded circle $\Bbb{S}^1$ in the Moore space:
\[
\Bbb{S}^1=\left\{ \left. e^{i\varphi }\right| 0\leq \varphi \leq \frac{2\pi }%
n\right\} \subset \Bbb{M}_n,
\]
while the quotient space is homeomorphic to the 2-sphere
$ \Bbb{M}_n/\Bbb{S}^1=\Bbb{S}^2$. The exact sequence of the pair
$\left( \Bbb{M}_n,\Bbb{S}^1\right) $ in $K$-theory reduces to
\[
\begin{array}{ccccc}
0\rightarrow  & K^1\left( \Bbb{S}^1\right)  &
\stackrel{\delta }{\longrightarrow } &
\widetilde{K}\left( \Bbb{S}^2\right)  & \rightarrow \widetilde{K}(\Bbb{M}_n)\to 0   \\
& \parallel  &  & \parallel  \\
  & \Bbb{ Z} &  & \Bbb{ Z}
\end{array}
\]
and the coboundary mapping
$K^1\left( \Bbb{S}^1\right) \stackrel{\delta }{\rightarrow }\widetilde{K}\left(
\Bbb{S}^2\right)$ acts as the multiplication by $n$. This description gives
\[
K^1\left( \Bbb{M}_n\right) =0,K\left( \Bbb{M}_n\right) =\Bbb{ Z\oplus Z}_n.
\]
The generator of the torsion part $\Bbb{ Z}_n$ is
$\left[ \gamma _n\right] -1\in K\left( \Bbb{M}_n\right)$, where
$\gamma _n\in \limfunc{Vect}\left( \Bbb{M}_n\right)$ is the pull-back of the
Hopf line bundle on ${\Bbb{S}^2}.$

The Whitney sum $n\gamma _n$ is a trivial vector bundle. In fact, its
transition function is equal to $\left( z,z,...,z\right) $. It is homotopic
to $\left( z^n,1,...,1\right) .$ Using the latter transition function, it
is easy to produce a trivialization
$n\gamma _n  \stackrel{\beta }{\longrightarrow }  \Bbb{ C}^n.$

{\bf 2. $K$-groups with coefficients.}
For a topological space $X$, the $K$-theory with coefficients $\Bbb{ Z}_n$ is
defined in terms of the usual (integral) $K$-theory by the formula
\begin{equation}
K^{*}\left( X,\Bbb{ Z}_n\right) =K^{*}\left( X\times \Bbb{M}_n,X\times pt\right).
\label{fdef}
\end{equation}
For instance, for a point we have
\[
K^*\left( pt,\Bbb{ Z}_n\right)=\widetilde{K}^*(\Bbb{M}_n),
\]
which is trivial for $K^1$ and $\Bbb{Z}_n$ for $K^1$.

There is an exact sequence in $K$-theory with coefficients
\begin{equation}
\begin{array}{c}
\rightarrow \!\!K\left( X\right) \!\stackrel{\times n}{\rightarrow }%
\!K\left( X\right) \!\rightarrow \!K\left( X,\Bbb{ Z}_n\right)
\!\!\rightarrow \!K^1\left( X\right) \!\stackrel{\times n}{\rightarrow }%
\!K^1\left( X\right) \!\rightarrow \!K^1\left( X,\Bbb{ Z}_n\right)
\!\rightarrow,
\end{array}
\label{sekq}
\end{equation}
corresponding to the short exact sequence
$0\to \Bbb{ Z}\stackrel{\times n}{\to }\Bbb{ Z\to Z}%
_n\to 0.$
It is obtained from the exact sequence of the pair $
\left( X\times \Bbb{M}_n,X\times \Bbb{S}^1\right) $ by the Bott periodicity.
We will need explicit descriptions of the connecting homomorphisms.
The ''reduction modulo $n$'' mappings
\[
K\left( X\right) \rightarrow K\left( X,\Bbb{ Z}_n\right) \text{\quad { %
and }\quad }K^1\left( X\right) \rightarrow K^1\left( X,\Bbb{ Z}_n\right)
\]
 are realized as tensor products with $\left[ \gamma _n\right]
-1\in \widetilde{K}\left( \Bbb{M}_n\right) ,${ \ while the Bokstein maps}
\[
K\left( X,\Bbb{ Z}_n\right) \rightarrow K^1\left( X\right) \text{\quad {
and }\quad }K^1\left( X,\Bbb{ Z}_n\right) \rightarrow K\left( X\right)
\]
{ are induced by the embedding } $\Bbb{S}^1\subset \Bbb{M}_n$.

\renewcommand{\refname}{References}
\addcontentsline{toc}{section}{References}


\vspace{1cm}

\hfill {\em Moscow, Potsdam}
\end{document}